\documentclass[12pt, reqno]{amsart}

\usepackage[utf8]{inputenc}
\usepackage[T2A]{fontenc}

\usepackage{amssymb,amscd,amsmath}
\usepackage[all]{xy}
\usepackage{graphicx}
\usepackage{color}
\usepackage{enumitem}
\usepackage{array}

\usepackage[margin=2.5cm]{geometry}

\newtheorem{theorem_eng}[subsection]{Theorem}
\newtheorem{lemma_eng}[subsection]{Lemma}

\newtheorem{proposition_eng}[subsection]{Proposition}
\newtheorem{corollary_eng}[subsection]{Corollary}

\newtheorem{definition_eng}[subsection]{Definition}

\newtheorem{remark_eng}[subsection]{Remark}

\makeatletter
\@addtoreset{subsection}{section}
\@addtoreset{equation}{section}
\@addtoreset{figure}{section}
\@addtoreset{table}{section}
\makeatother

\makeatletter
\newcommand\testshape{family=\f@family; series=\f@series; shape=\f@shape.}
\def\myemphInternal#1{\if n\f@shape%
\begingroup\itshape #1\endgroup\/%
\else\begingroup\bfseries #1\endgroup%
\fi}
\def\myemph{\futurelet\testchar\MaybeOptArgmyemph}
\def\MaybeOptArgmyemph{\ifx[\testchar \let\next\OptArgmyemph
                 \else \let\next\NoOptArgmyemph \fi \next}
\def\OptArgmyemph[#1]#2{\index{#1}\myemphInternal{#2}}
\def\NoOptArgmyemph#1{\myemphInternal{#1}}
\makeatother

\newcommand{\bN}{\mathbb{N}}
\newcommand{\bZ}{\mathbb{Z}}
\newcommand{\bR}{\mathbb{R}}

\newcommand{\Int}{\mathop{\mathrm{Int}}\nolimits}

\newcommand{\Cl}[1]{\overline{#1}}

\renewcommand{\emptyset}{\varnothing}

\newcommand\Xsp{X}
\newcommand\Ysp{Y}
\newcommand\Zsp{Z}
\newcommand\Usp{U}
\newcommand\Vsp{V}
\newcommand\Wsp{W}
\newcommand\Partition{\Delta}

\newcommand\prj{p}
\newcommand\bnd[1]{\mathrm{hcl}(#1)}

\newcommand\YspecPtSet{\mathcal{V}}

\newcommand\eps{\varepsilon}

\newcommand{\crs}{\gamma}  
\newcommand{\crsSubDomain}{\Wsp}  
\newcommand{\crsDomain}{\Vsp}  

\newcommand\HalfSpace[1]{\bR^{#1}_{+}}

\newcommand\Kset[2]{K(#1, #2)}
\newcommand\Lset[2]{L(#1, #2)}

\newcommand{\pa}{u}
\newcommand{\pb}{v}
\newcommand{\pc}{w}
\newcommand{\ps}{s}
\newcommand{\pt}{t}

\newcommand{\px}{x}
\newcommand{\py}{y}

\newcommand\leaf{\mathbf{\omega}}
\newcommand\spleaf{\sigma}
\newcommand\pleaf[1]{\leaf_{{#1}}}

\newcommand\segm[1]{I_{#1}}

\newcommand\chartMap{\varphi}
\newcommand\folEmb{\psi}
\newcommand\bijleaf{\phi}

\newcommand\SAT[1]{\mathbf{S}(#1)}

\newcommand\Pel{P}
\newcommand\Qel{Q}
\newcommand\Family{\mathcal{Z}}

\newcommand\acrs{\mathbf{A}}
\newcommand\bcrs{\mathbf{B}}

\newcommand{\Btrans}{B^{n}}

\newcommand\dblX{\widehat{\Xsp}}
\newcommand\dblV{\widehat{\crsDomain}}
\newcommand\dblP{\widehat{\Partition}}
\newcommand\dblc{\widehat{\crs}}

\author{Sergiy Maksymenko, Eugene Polulyakh}
\title{One-dimensional foliations on topological manifolds}

\address{Topology Laboratory, Institute of Mathematics, Ukrainian National Academy of Science, Te\-re\-shchen\-kiv\-ska str. 3, 01004 Kyiv, Ukraine}
\email{maks@imath.kiev.ua, polulyah@imath.kiev.ua}



\begin{document}

\begin{abstract}
Let $X$ be an $(n+1)$-dimensional manifold, $\Delta$ be a one-dimensional foliation on $X$, and $p: X \to X / \Delta$ be a quotient map.
We will say that a leaf $\omega$ of $\Delta$ is \textit{special} whenever the space of leaves $X / \Delta$ is not Hausdorff at $\omega$.
We present necessary and sufficient conditions for the map $p: X \to X / \Delta$ to be a locally trivial fibration under assumptions that all leaves of $\Delta$ are non-compact and the family of all special leaves of $\Delta$ is locally finite.
\end{abstract}

\maketitle

\section{Introduction}
Study of the topological structure of flow lines foliations has a long history and leads back to H.~Poincar\`e.
The question when a partition into curves is a foliation was considered by H.~Whitney~\cite{Whitney:AM:1933}, \cite{Whitney:BAMS:1941}.
In two-dimensional case one-dimensional foliations appeared as level-sets of pseudo-harmonic functions in W.~Kaplan~\cite{Kaplan:DJM:1940}, \cite{Kaplan:DJM:1941}.

Let $\Partition$ be a one-dimensional foliation on $\bR^2$, and $\bR^2/\Partition$ be the space of leaves endowed with the quotient topology.
Notice that $\bR^2/\Partition$ is usually non-Hausdorff.
W.~Kaplan~\cite{Kaplan:DJM:1940} showed that
\begin{enumerate}[leftmargin=*, label=$(\arabic*)$]
\item\label{enum:Kaplan:p_is_loc_triv}
the quotient map $\prj:\bR^2\to\bR^2/\Partition$ is a locally trivial fibration with fiber $\bR$;
\item\label{enum:Kaplan:split_into_strips}
there exists at most countably many leaves $\{\leaf_i\}_{i\in A}$ of $\Partition$ such that the complement $\bR^2\setminus\{\leaf_i\}_{i\in A}$ is a disjoint union $\mathop{\sqcup}\limits_{j\in B} S_j$, where each $S_j$ is homeomorphic with $(0,1)\times \bR$ so that the lines $t\times\bR$, $t\in(0,1)$, correspond to the leaves of $\Partition$;
\item\label{enum:Kaplan:pseudoharm_func}
there exists a pseudoharmonic function (without singularities) $f:\bR^2\to\bR$ whose foliation by connected components of level-sets coincides with $\Partition$.
\end{enumerate}
See also W.~Boothby~\cite{Boothby:AJM_1:1951}, \cite{Boothby:AJM_2:1951}, M.~Morse and J.~Jenkins~\cite{JenkinsMorse:AJM:1952}, \cite{JenkinsMorse:APNASUSA:1953}, \cite{JenkinsMorse:ActaMath:1954}, \cite{JenkinsMorse:UMichPress:1955} and M.~Morse~\cite{Morse:FM:1952}, \cite{Morse:JMPA:1956} for extensions of Kaplan's results to foliations with singularities.

A.~Haefliger and G.~Reeb~\cite{HaefligerReeb:EM:1957} studied general one-dimensional \textit{non-Hausdorff} manifolds and showed, in particular, that the above result~\ref{enum:Kaplan:pseudoharm_func} of W.~Kaplan can be deduced from Poincar\`e-Bendixon theorem, see also~\cite{Haefliger:CTS:1954}, \cite{Reeb:CTS:1954}.

Later C.~Godbillon and G.~Reeb~\cite{GodbillonReeb:EM:1966} classified locally trivial fibrations over a \myemph{non-Hausdorff letter $Y$}.
Though they considered a very special case their methods clarify the general situation.

The question when for an arbitrary $k$-dimensional foliation $\Partition$ on $\Xsp$ the quotient map $\prj:\Xsp\to\Xsp/\Partition$ has homotopy lifting properties was considered in C.~Godbillon~\cite{Godbillon:AIF:1967}, see also G.~Meigniez~\cite{Meigniez:AIF:1997} and \cite{Meigniez:TrAMS:2002} for the criterion when $\prj$ is a Serre fibration or a locally trivial fibration but mostly in smooth category.
J.~Harrison~\cite{Harrison:Top:1988} studied similar problem concerning geodesic flows without compact orbits.
Also foliations by flow lines on $3$-manifolds are classified by S.~Matsumoto~\cite{Matsumoto:MAMS:2003}.

In recent years a progress in the theory of Hamiltonial dynamical systems of small degrees of freedom increased an interest to the structure of level-sets functions on surfaces, see e.g. A.~Fomenko and A.~Bolsinov~\cite{BolsinovFomenko:1997}, A.~Oshemkov~\cite{Oshemkov:PSIM:1995}, V.~Sharko~\cite{Sharko:UMZ:2003}, \cite{Sharko:Zb:2006}, V.~Sharko and Yu.~Soroka~\cite{SharkoSoroka:MFAT:2015}, E.~Polulyakh and I.~Yurchuk~\cite{PolulyakhYurchuk:ProcIM:2009}, E.~Polulyakh~\cite{Polulyakh:UMZ:ENG:2015}.

Homotopy properties of foliations on surfaces glued from strips similarly to~\ref{enum:Kaplan:split_into_strips} are studied in S.~Maksymenko and E.~Polulyakh~\cite{MaksymenkoPolulyakh:PGC:2015} and \cite{MaksymenkoPolulyakh:MFAT:2016} and Yu.~Soroka~\cite{Soroka:MFAT:2016}.

In \cite{MaksymenkoPolulyakh:MFAT:2016} the authors extended Kaplan's result~\ref{enum:Kaplan:split_into_strips} to foliations on arbitrary non-compact surfaces $\Xsp$.
Namely, under certain assumptions including~\ref{enum:Kaplan:p_is_loc_triv}, i.e. that $\prj:\Xsp\to\Xsp/\Partition$ is a locally trivial fibration, the topological structure of the closures $\overline{S_j}$ of strips $S_j$ was described.

In the present paper we consider an arbitrary one-dimensional foliation $\Partition$ with all non-compact leaves on a topological manifold $\Xsp$.
Our main result gives necessary and sufficient conditions for the quotient map $\prj:\Xsp\to\Xsp/\Partition$ to be a locally trivial fibration, see Theorem~\ref{th:charact_loc_triv}.

As mentioned above such types of questions were extensively studied.
However, the essentially new features of Theorem~\ref{th:charact_loc_triv} in comparison e.g.\! with \cite{Godbillon:AIF:1967}, \cite{Meigniez:TrAMS:2002} and others, is that we work in $C^0$ category only and give a characterization in terms of the topology of the quotient space $\Xsp/\Partition$.

\section{One-dimensional foliations}\label{sect:one_dim_foliations}
Let $\HalfSpace{n} = \{ (x_1,\ldots,x_n) \mid x_n \geq0 \}$ be the closed half-space in $\bR^{n+1}$. 
\begin{definition_eng}[cf.~\cite{CandelConlon:AMS:2000}]\label{defn:foliated_chart}
Let $\Xsp$ be an $(n+1)$-dimensional topological manifold, $n \geq 1$.
A \myemph{foliated chart} on $X$ of codimension $n$ is a pair $(\Usp, \chartMap)$, where $\Usp \subset \Xsp$ is open and $\chartMap : \Usp \to (a,b) \times \Btrans$ is a homeomorphism with $\Btrans$ being an open subset of $\HalfSpace{n}$ and $a<b\in\bR\cup\{\pm\infty\}$.
The set $P_y = \chartMap^{-1}\bigl((a,b) \times \{y\}\bigr)$, $y \in \Btrans$, is called a \myemph{plaque} of this foliated chart.
\end{definition_eng}

\begin{definition_eng}[cf.~\cite{CandelConlon:AMS:2000,Tamura:1979}]\label{defn:1-foliation}
Let $\Partition = \{ \leaf_\alpha \,|\, \alpha \in A \}$ be a partition of $X$ into pathconnected subsets $\leaf_\alpha$ of $X$.
Suppose that $X$ admits an atlas $\{\Usp_i, \chartMap_i\}_{i\in\Lambda}$ of foliated charts of codimension $n$ such that, for each $\alpha \in A$ and each $i \in \Lambda$, every pathcomponent of a set $\leaf_\alpha \cap \Usp_i$ is a plaque.
Then $\Partition$ is said to be a \myemph{foliation of $X$ of dimension 1} (and \myemph{codimension $n$}) and $\{\Usp_i, \chartMap_i\}_{i\in\Lambda}$ is called a \myemph{foliated atlas associated to $\Partition$}.
Each $\leaf_\alpha$ is called a \myemph{leaf} of the foliation and the pair $(X, \Partition)$ is called a \myemph{foliated manifold}.
\end{definition_eng}

\begin{remark_eng}\rm
In~\cite{Kaplan:DJM:1940} one-dimensional foliations on the plane were also called \myemph{regular families of curves}.
\end{remark_eng}

In what follows we will assume that $\Xsp$ is endowed with some $1$-dimensional foliation $\Partition$.
We will also consider only foliated charts included into some (maximal) foliated atlas associated to $\Partition$.

Let $\leaf$ be a leaf of $\Partition$.
If $\leaf$ is compact, then it is homeomorphic with the circle.
Otherwise, there exists a continuous bijection $\bijleaf:\bR\to\leaf$.
Moreover, if $\bar{\bijleaf}:\bR\to\leaf$ is another continuous bijection, then $\bijleaf^{-1}\circ\bar{\bijleaf}:\bR\to\bR$ is a homeomorphism, c.f.~\cite[Proposition~6]{Godbillon:AIF:1967}.

Recall that a continuous map $f:A\to B$ is called \myemph{proper} whenever for each compact $K\subset B$ its inverse image $f^{-1}(K)$ is compact.
The following lemma is easy and we leave it for the reader.
\begin{lemma_eng}\label{lm:proper_leaves}
Consider the following conditions on $\leaf \in \Partition$:
\begin{itemize}[leftmargin=3em]
\item[$\mathrm{(m)}$]
there exists an embedding $\bijleaf:\bR\to\Xsp$ with $\bijleaf(\bR) = \omega$;
\item[$\mathrm{(p)}$]
there exists a proper injective continuous map $\bijleaf:\bR\to\Xsp$ with $\bijleaf(\bR) = \omega$;
\item[$\mathrm{(p)'}$]
any injective continuous map $\bijleaf:\bR\to\Xsp$ with $\bijleaf(\bR) = \omega$ is proper;
\item[$\mathrm{(c)}$]
$\omega$ is a closed subset of $\Xsp$.
\end{itemize}
Then the following equivalences hold true:
\[ \mathrm{(m)} \ \& \ \mathrm{(c)} \ \Leftrightarrow \ \mathrm{(p)} \ \Leftrightarrow \ \mathrm{(p')}.\]
\end{lemma_eng}

A leaf $\leaf$ satisfying condition (p) of Lemma~\ref{lm:proper_leaves} will be said to be \myemph{properly embedded}%
\footnote{In the book~\cite[\S16]{Tamura:1979} a leaf is called \myemph{proper} if it satisfies condition (m).}.

The union of all leaves of $\Partition$ intersecting a subset $\Usp\subset\Xsp$ is called the \myemph{saturation} of $\Usp$ and denoted by $\SAT{\Usp}$.
The following lemma is easy to prove.
\begin{lemma_eng}\label{lm:saturation_is_open}{\rm\cite[Proposition~1.5]{Godbillon:F:1991}, \cite[Theorem~4.10]{Tamura:1979}}
If $\Usp \subset \Xsp$ is open, then $\SAT{\Usp}$ is open as well.
\qed
\end{lemma_eng}

If $\Ysp$ is a manifold, then a \myemph{trivial} $1$-dimensional foliation $\Partition$ on the product $\bR\times\Ysp$ is a partition of $\bR\times\Ysp$ into the lines $\bR\times y$, $y\in\Ysp$.

Let $\Partition_i$ be a $1$-foliation on $\Xsp_i$, $i=1,2$.
Then an embedding $\folEmb:\Xsp_1 \to \Xsp_2$ will be called \myemph{foliated} whenever $\folEmb(\leaf)$ is contained in some leaf of $\Partition_2$ for each leaf $\leaf\in\Partition_1$.

In particular, if $\chartMap : \Usp \to (a,b) \times \Btrans$ is a foliated chart as in Definition~\ref{defn:foliated_chart}, then its inverse $\folEmb = \chartMap^{-1}:(a,b) \times \Btrans \to \Xsp$ is an open foliated embedding.
In this case the set $P_{\pa} = \folEmb\bigl( (a,b) \times \{\pa\}\bigr)$ is a plaque for each $\pa\in\Btrans$.

\medskip 
\noindent
{\bf Space of leaves.}
Let $\Ysp = \Xsp/\Partition$ be the space of leaves and $\prj:\Xsp\to\Ysp$ be the corresponding quotient map.
Endow $\Ysp$ with the \emph{quotient topology} with respect to $\prj$.
Thus a subset $\Vsp\subset\Ysp$ is open if and only if its inverse $\prj^{-1}(\Vsp)$ is open in $\Xsp$.

Notice that for a subset $\Usp \subset \Xsp$ its saturation is $\SAT{\Usp} = \prj^{-1}\bigl(\prj(\Usp)\bigr)$.
In particular, Lemma~\ref{lm:saturation_is_open} means that $\prj$ is an open map.

Evidently, $\Ysp$ is a $T_1$-space if and only if each leaf of $\Partition$ is a closed subset of $\Xsp$.
However, in general, $\Ysp$ is not a Hausdorff space.

\medskip
\noindent
{\bf Special points.}
Let $\pa\in\Ysp$ be a point and $\beta_{\pa}$ be a base of neighborhoods of $\pa$.
Then the following set
\[
\bnd{\pa} := \mathop{\cap}\limits_{\Vsp\in\beta_{\pa}} \overline{\Vsp}
\]
will be called the \myemph{Hausdorff closure} of $\pa$.
A point $\pa$ will be called \myemph{special}\footnote{In~\cite[Definition~3]{HaefligerReeb:EM:1957} such a point is called a \myemph{branch} point.
See also~\cite{GodbillonReeb:EM:1966}.} if $\pa\not=\bnd{\pa}$.
Notice that $\pa\in\bnd{\pb}$ if and only if any two neighborhoods of $\pa$ and $\pb$ intersect.
The latter statement is symmetric with respect to $\pa$ and $\pb$, and so it is equivalent to the assumption $\pb\in\bnd{\pa}$.
However, one easily checks that the property ``belong to Hausdorff closure'' is not transitive.

Evidently, $\Ysp$ is Hausdorff if and only if $\pa=\bnd{\pa}$ for all $\pa\in\Ysp$, that is when $\Ysp$ has no special points.
The set of all special points of $\Ysp$ will be denoted by $\YspecPtSet$.

We will say that a leaf $\leaf$ of $\Partition$ is \myemph{special} if $\prj(\leaf)$ is a special point of $\Ysp$.
In particular, $\Sigma :=\prj^{-1}(\YspecPtSet)$ is the set of all special leaves of $\Partition$.

The following lemma gives a characterization of special leaves and extends~\cite[Proposition~4]{GodbillonReeb:EM:1966}.
\begin{lemma_eng}\label{lm:charact_spec_leaves}
Let $\leaf\in\Partition$ be a leaf and $\pa = \prj(\leaf)$ be the corresponding point in $\Ysp$.
Then the following conditions equivalent:
\begin{enumerate}[leftmargin=2em]
\item[\rm(1)]
$\pa$ is a special point of $\Ysp$, and so $\leaf$ is a special leaf of $\Partition$;

\item[\rm(2)]
there exists a point $\pb\in\bnd{\pa}$ distinct from $\pa$ and a sequence $\{\pc_{i}\}_{i\in\bN}$ converging to both points $\pa$ and $\pb$;

\item[\rm(3)]
there exist two sequences $\{\px_{i}\}_{i\in\bN}$ and $\{\py_{i}\}_{i\in\bN}$ in $\Xsp$ such that $\px_{i}$ and $\py_{i}$ belong to the same leaf for all $i\in\bN$, that is $\prj(\px_{i})=\prj(\py_{i})$, and
\begin{align}\label{equ:spec_sequences}
 \lim\limits_{i\to\infty} \px_{i} &= \px \in \leaf, &
 \lim\limits_{i\to\infty} \py_{i} &= \py \not\in \leaf.
\end{align}
\end{enumerate}
\end{lemma_eng}
\begin{proof}
Equivalence (1)$\Leftrightarrow$(2) is well known and easy.

(3)$\Rightarrow$(2).
Denote $w_i = \prj(\px_{i})=\prj(\py_{i})$, $i\in\bN$ and $\pb = \prj(\py)$.
Then, by continuity of $\prj$, the sequence $\{\pc_{i}\}_{i\in\bN}$ converges to distinct points $\pa$ and $\pb$.
In particular, $\pb\in\bnd{\pa}$.

(2)$\Rightarrow$(3).
Choose any points $\px\in\leaf$ and $\py\in\leaf':=\prj^{-1}(\pb)$, and let $\{\Usp_i\}_{i\in\bN}$ and $\{\Vsp_i\}_{i\in\bN}$ be countable bases of topology on $\Xsp$ at $\px$ and $\py$ respectively.
Since $\prj$ is open, $\prj(\Usp_i)$ and $\prj(\Vsp_i)$ are open neighborhoods of $\pa$ and $\pb$ respectively.
But these points are special and $\pa\in\bnd{\pb}$, whence $\prj(\Usp_i) \cap \prj(\Vsp_i) \not=\varnothing$ for all $i\in\bN$.
Hence there exist $\px_{i} \in \Usp_i$ and $\py_{i}\in\Vsp_i$ such that $\prj(\px_{i}) = \prj(\py_{i})$.
Then $\{\px_{i}\}_{i\in\bN}$ and $\{\py_{i}\}_{i\in\bN}$ converge to $\px$ and $\py$ respectively.
\end{proof}

\begin{definition_eng}\label{def:cross_section}
Let $\crsDomain$ be an open subset of $\bR^n_{+}$, where $n=\dim\Xsp-1$ is the codimension of the foliation $\Partition$.
A continuous map $\crs:\crsDomain\to\Xsp$ is called a \myemph{cross section} of $\Partition$ if
\begin{itemize}[leftmargin=1.1em, label=$\bullet$]
\item $\crs(\crsDomain \cap \partial\bR^n_{+}) \subset \partial\Xsp$;
\item the composition $\prj\circ\crs:\crsDomain\to\Ysp$ is an injective map, that is $\pa\not=\pb\in\crsDomain$ implies that the images $\crs(\pa)$ and $\crs(\pb)$ of these points belong to distinct leaves of $\Partition$.
\end{itemize}
If $\px\in\crs(\crsDomain)$ and $\leaf$ is a leaf of $\Partition$ containing $\px$, then we will also say that $\crs$ \myemph{passes through $\px$} as well as \myemph{through $\leaf$}.
\end{definition_eng}

The aim of this paper is to present necessary and sufficient conditions for the map $\prj$ to be a locally trivial fibration under assumption that all leaves of $\Partition$ are non-compact.

\begin{theorem_eng}\label{th:charact_loc_triv}
Let $\Xsp$ be an $(n+1)$-dimensional manifold and $\Partition$ be a one-dimensional foliation on $\Xsp$.
Suppose that all leaves of $\Partition$ are non-compact and the family of all special leaves of $\Partition$ is locally finite.
Then the following conditions are equivalent.
\begin{enumerate}[leftmargin=*, label=$(\arabic*)$]
\item\label{lm:loc_triv_relations:enum:prj_is_ltfibr}
The quotient map $\prj:\Xsp\to\Xsp/\Partition$ is a locally trivial fibration with fiber $\bR$ and $\Ysp$ is locally homeomorphic with $\bR^n_{+}$ (though it is not necessary a Hausdorff space).
\item\label{lm:loc_triv_relations:enum:partit_lt__prj_open}
For each leaf $\leaf$ there exists an open saturated neighborhood foliated homeomorphic with $\bR\times \Vsp$, where $\Vsp$ is an open subset of $\bR^n_{+}$.
\item\label{lm:loc_triv_relations:leaf_proper_cross_sections}
For each leaf $\leaf$ of $\Partition$ there exists a cross section passing through $\leaf$.
\end{enumerate}
\end{theorem_eng}

\begin{remark_eng}\rm
It is proved in~\cite[Chapter III, Propostion~4 and Corollary]{Godbillon:AIF:1967} that for arbitrary $p$-dimensional foliation $\Partition$ then the quotient map $\prj:\Xsp\to\Xsp/\Partition$ is a Serre fibration whenever it satisfies a certain variant of homotopy extension property and either has a local section at each point or the quotient $\Xsp/\Partition$ is a (possibly non-Hausdorff manifold).
See also~\cite{Meigniez:AIF:1997} and~\cite{Meigniez:TrAMS:2002} for extensions.

Our Theorem~\ref{th:charact_loc_triv} claims that for one-dimensional foliations $\Partition$  with locally finite family of special leaves existence of cross sections with open subsets of $\bR_{+}^{n}$ implies that $\prj$ is even a locally trivial fibration and $\Xsp/\Partition$ is a possibly non-Hausdorff manifold.
\end{remark_eng}

\begin{remark_eng}\rm
Equivalence between~\ref{lm:loc_triv_relations:enum:prj_is_ltfibr} and~\ref{lm:loc_triv_relations:enum:partit_lt__prj_open} for $\dim\Xsp=2$ is proved in~\cite{MaksymenkoPolulyakh:MFAT:2016} without assumption that $\Ysp$ is locally homeomorphic with $\bR$.
Also in~\cite[\S2.2, Proposition~1]{HaefligerReeb:EM:1957} it is show that $\Xsp/\Partition$ is a $1$-manifold for one-dimensional foliation on $\bR^2$.
\end{remark_eng}

\begin{remark_eng}\rm
R.~H.~Bing~\cite{Bing:BAMS:1958}, \cite{Bing:AM:1959} constructed a non-manifold $B \subset \bR^4$ such that $\bR\times B$ is homeomorphic with $\bR^4$.
In other words, $\bR^4$ admits a trivial partition into open arcs (being not a foliation) such that the quotient space $B$ is not a $3$-manifold.
That example was improved by many authors, see e.g R.~Rosen~\cite{Rosen:AM:1961}, J.~Kim~\cite{Kim:ProcAMS:1964}, J.~Bailey~\cite{Bailey:TrAMS:1970}, L.~Rubin~\cite{Rubin:TrAMS:1972}.
\end{remark_eng}

\begin{remark_eng}\rm
E.~Dyer and M.~Hamstrom~\cite{DyerHamstrom:FM:1958} studied so called \myemph{completely regular} mappings $\prj:\Xsp\to\Ysp$ between metric spaces such that the inverse images of all points are in a certain sense ``uniformly homeomorphic'', and get sufficient conditions when such a map is equivalent to a trivial fibration, see~\cite[Theorem~7]{DyerHamstrom:FM:1958}, and also \cite{McAuley:TSW:1966}, \cite{Seidman:TrAMS:1970} for generalizations.
We consider here a similar problem, but now the space $\Ysp$ is not even Hausdorff, and we gave conditions when $\prj$ is a locally trivial fibration.
\end{remark_eng}

The following statement is proved in~\cite[Theorem~1]{SharkoSoroka:MFAT:2015} for continuous functions $f:\bR^2\to\bR$, and in~\cite[item 3 at the end of page 3778]{Meigniez:TrAMS:2002} for smooth case.
\begin{theorem_eng}
Let $M$ and $N$ be two manifolds such that $\dim M = \dim N + 1$ and $f:M \to N$ be a surjective continuous map such that
\begin{itemize}
\item
$f(\Int{M}) = \Int{N}$ and $f(\partial M) = \partial N$;
\item
the partition $\Partition = \{ f^{-1}(c) \mid c\in N\}$ of $M$ constitutes a one-dimensional foliation with all non-compact leaves.
\end{itemize}
Then $f$ is a locally trivial fibration with fiber $\bR$.
In particular, if $N$ is contractible, then $f$ is a trivial fibration.
\end{theorem_eng}
\begin{proof}
We claim that $\Partition$ contains no special leaves and each leaf admits a cross section.
Then it will follow from Theorem~\ref{th:charact_loc_triv} that $f$ is a locally trivial fibration with fiber $\bR$.

{\em Absence of special leaves.}
Let $Y = M / \Partition$ be the space of leaves endowed with the corresponding factor topology.
Then $f$ can be written as a composition of the following maps
\[f = \theta \circ \prj: M \xrightarrow{~~\prj~~} Y \xrightarrow{~~\theta~~} N,\]
where $\theta$ is the induced continuous bijection.
Since $N$ is Hausdorff, it follows that so is $Y$, and therefore $Y$ contains no special points.
Hence $\Partition$ contains no special leaves.

{\em Existence of cross sections.}
Let $\px\in M$ and $\chartMap: \Usp \to (-1,1)\times \Btrans$ be a foliated chart at $\px$ as in Definition~\ref{defn:1-foliation} such that $\chartMap(\px) = (0,0) \in (-1,1)\times \Btrans$, where $n=\dim N$.
Then the map $\crs:\Btrans \to M$ defined by $\crs(\py) = \chartMap^{-1}(0,\py)$ is a cross section of $\Partition$.
\end{proof}

In fact, Theorem~\ref{th:charact_loc_triv} is an easy consequence of the following statements:

\noindent{\bf Lemma~\ref{lemma:proper_leaf1}.}~{\em 
Let $\leaf_0$ be a leaf of $\Partition$. Suppose that for each leaf $\leaf$ of $\Partition$ contained in $\SAT{\Cl{\leaf_0}}$ there exists a cross section $\crs$ passing through $\leaf$.
Then $\leaf_0$ is properly embedded.
}

\begin{theorem_eng}\label{th:trivialization_of_cross_section}
Let $\crs:\crsDomain \to \Xsp$ be a cross section intersecting only leaves being simultaneously non-compact, properly embedded, and non-special.
Then the saturation $\SAT{\crs(\crsDomain)}$ is open and foliated homeomorphic with $\bR\times\crsDomain$.
\end{theorem_eng}

The proof of Theorem~\ref{th:charact_loc_triv} will be given in \S\ref{sect:proof_of_theorem:th:charact_loc_triv}.
In~\S\ref{sect:prelim} we will prove some general preliminary results concerning one-dimensional $C^0$ foliations being well known for smooth case.
In particular we will prove Lemma~\ref{lemma:proper_leaf1}.
\S\ref{sect:proof_of_theorem:th:trivialization_of_cross_section} is devoted to the proof of Theorem~\ref{th:trivialization_of_cross_section} using E.~Michael's theorems about selections of multivalued maps.

\section{Proof of Theorem~\ref{th:charact_loc_triv}}\label{sect:proof_of_theorem:th:charact_loc_triv}
\ref{lm:loc_triv_relations:enum:prj_is_ltfibr}$\Rightarrow$\ref{lm:loc_triv_relations:enum:partit_lt__prj_open}, \ref{lm:loc_triv_relations:leaf_proper_cross_sections}.
Suppose the quotient map $\prj:\Xsp\to\Ysp$ is a locally trivial fibration with fiber $\bR$ and $\Ysp$ is locally homeomorphic with $\bR^n_{+}$.
This means that for each $\leaf\in\Partition$ there exist
\begin{itemize}
\item an open neighborhood $\Vsp \subset \Ysp$ of its image $\pa=\prj(\leaf)$ homeomorphic with an open subset of $\bR^n_{+}$ and
\item a foliated homeomorphism $\folEmb:\bR\times\Vsp\to \prj^{-1}(\Vsp)$.
\end{itemize}
Then $\prj^{-1}(\Vsp)$ is an open and saturated neighborhood of $\leaf$ and $\folEmb$ is a foliated homeomorphism required by~\ref{lm:loc_triv_relations:enum:partit_lt__prj_open}.

Moreover, the map $\crs:\Vsp \to \Xsp$ defined by $\crs(\pb) = \folEmb(0,\pb)$ is a cross section passing through $\leaf$.
This proves \ref{lm:loc_triv_relations:leaf_proper_cross_sections}.

\smallskip

\ref{lm:loc_triv_relations:leaf_proper_cross_sections}$\Rightarrow$\ref{lm:loc_triv_relations:enum:partit_lt__prj_open}.
Suppose each leaf of $\Partition$ admits a local cross section.
Then it follows from Lemma~\ref{lemma:proper_leaf1} that all leaves of $\Partition$ are properly embedded.
Let $\Sigma$ be a family of all special leaves and $\spleaf\in\Sigma$ be a special leaf.

Since each leaf is closed and $\Sigma$ is a locally finite family, it follows that $\Sigma\setminus\spleaf$ is a closed set, whence $\Xsp' = (\Xsp\setminus\Sigma)\cup\spleaf$ is open and saturated and contains no special leaves.
Moreover, since each leaf in $\Xsp'$ admits a local cross section, it follows from Theorem~\ref{th:trivialization_of_cross_section} that each leaf $\leaf \subset \Xsp'$ has an open saturated neighborhood $\Wsp$ foliated homeomorphic with $\bR\times\Vsp$, where $\Vsp$ is an open subset of $\bR^{n}_{+}$.
Then $\Wsp$ is also open in $\Xsp$.
This proves~\ref{lm:loc_triv_relations:enum:partit_lt__prj_open}.

\smallskip

\ref{lm:loc_triv_relations:enum:partit_lt__prj_open}$\Rightarrow$\ref{lm:loc_triv_relations:enum:prj_is_ltfibr}.
Let $\pa\in\Ysp$ and $\leaf = \prj^{-1}(\pa)$ be the corresponding leaf of $\Partition$.
Suppose there exist an open $\Vsp \subset \bR^{n}_{+}$ and a foliated homeomorphism $\folEmb:\bR\times\Vsp\to \Wsp_{\leaf}$ onto some open and saturated neighborhood $\Wsp_{\leaf}$ of $\leaf$.
Since $\prj$ is an open map, so is the composition $\prj\circ\folEmb$.
Hence $\Usp_{\pa} := \prj\bigl(\folEmb(\Wsp_{\leaf})\bigr)$ is an open neighborhood of $\pa$ in $\Ysp$.
Moreover, the restriction $\prj\circ\folEmb|_{0\times \bR^n}: 0\times \Vsp \to \Usp_{\pa}$ is a continuous and open bijection, and so it is a homeomorphism.
Thus $\Ysp$ is locally homeomorphic with $\bR^{n}_{+}$ and the map $\prj\circ\folEmb:\bR\times\Vsp\to\Usp_{\pa}$ is a trivialization of $\prj$ over $\Usp_{\pa}$, so $\prj$ is a locally trivial fibration with fiber $\bR$.
Theorem~\ref{th:charact_loc_triv} is completed.

\section{Preliminaries}\label{sect:prelim}

In this section we will assume that $\Xsp$ is an $(n+1)$-dimensional manifold with $\partial\Xsp=\varnothing$ and $\Partition$ is a one-dimensional foliation on $\Xsp$.

Some statements in this section are well known for $C^1$ foliations e.g.~\cite{Meigniez:TrAMS:2002}, and some of them are proved for $C^0$ case but for the foliations on $\bR^2$, see e.g. W.~Kaplan~\cite{Kaplan:DJM:1940}.
However we did not find good exposition in the literature for general $C^0$ foliations needed in our case and therefore short proofs will be presented.
This will also make the paper self-contained.

It will be convenient to regard the graph of a function $f:\Xsp\to\bR$ as the following subset \[\Gamma_f:= \{(f(x),x) \mid x\in\Xsp\}\] of $\bR\times\Xsp$.
Thus we switch the coordinates.

\begin{figure}[h]
\includegraphics[height=1.2cm]{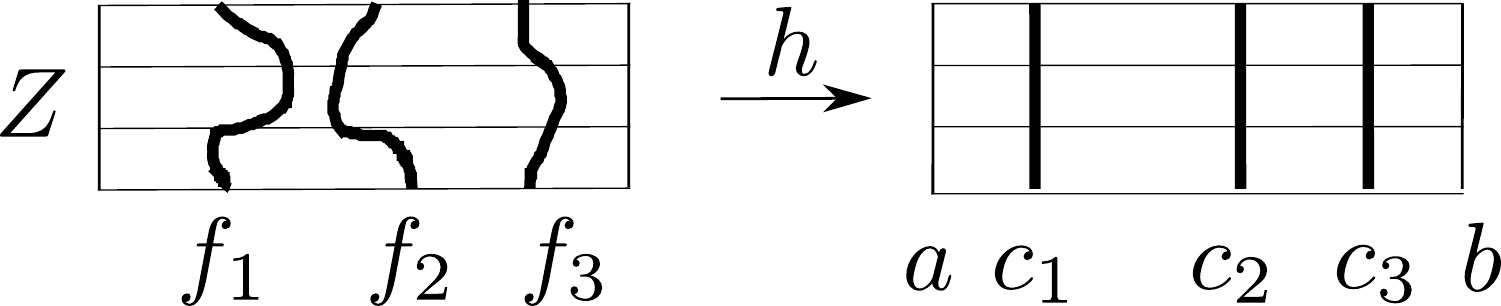}
\caption{}\label{fig:const_graphs}
\end{figure}

\begin{lemma_eng}\label{lm:make_constrant_graph}
Let $\Zsp$ be a topological space, $f_1,\ldots,f_k:\Zsp\to(a,b)$ be continuous functions such that $f_i(z) < f_j(z)$ for all $i<j$ and $z\in\Zsp$, and
\[ \Gamma_i = \{ (f_i(z), z) \mid z\in\Zsp \} \ \subset \ [a,b]\ \times \Zsp \]
be the graph of $f_i$.
Let also $c_1<c_2<\cdots < c_k \in (a,b)$ be any increasing $k$-tuple of numbers.
Then there exists a self-homeomorphism $h$ of $[a,b] \times \Zsp$ such that
\begin{enumerate}[leftmargin=*, label=$\mathrm{(\alph*)}$]
\item
$h$ is fixed on $a \times  \Zsp$ and $b \times \Zsp$;
\item\label{lm:make_constrant_graph:enum:preserve_lines}
$h([a,b]\times z) = [a,b] \times z$ for all $z\in\Zsp$;
\item
$h(\Gamma_i) = c_i \times \Zsp$;
\item\label{lm:make_constrant_graph:enum:fixed_on_0}
if $f_i(z) = c_i$ for some $z\in \Zsp$ and all $i=1,\ldots,k$, then $h$ is fixed on $[a,b]\times z$.
\end{enumerate}
\end{lemma_eng}
\begin{proof}
The proof follows from~\cite[Lemma~6.1.1]{Maksymenko:BSM:2006}, see also~\cite[Lemma 5.2.1]{MaksymenkoPolulyakh:MFAT:2016}.
Let us just mention that the situation can be reduced to the case $[a,b] = [0,1]$, and that for $k=1$ the desired self-homeomorphism $h$ of $[0,1] \times \Zsp$ can be defined e.g.\! by
\[
h(s,z) =
\begin{cases}
(s, z), & s\in\{0,1\}, \\
\bigl(s^{\log_{f_1(z)}c_1}, z\bigr), & s\in(0,1).
\end{cases}
\]
We leave the details for the reader, see Figure~\ref{fig:const_graphs}.
\end{proof}

\begin{lemma_eng}\label{lm:adopt_quasi_tube}
Let $W$ be an open neighborhood of $0$ in $\bR^n$ and $\crs:W \to (a,b)\times \bR^n$ be a cross section of the trivial one-dimensional foliation such that $\crs(0) \in (a,b)\times 0$.
Then for each $c\in (a,b)$ there exists an open embedding $\folEmb:(a,b)\times W \subset (a,b)\times \bR^n$ such that
\begin{enumerate}[leftmargin=*, label=$\mathrm{(\roman*)}$]
\item\label{lm:adopt_quasi_tube:enum:pres_lines}
$\folEmb\bigl( (a,b)\times x\bigr) = (a,b)\times\crs(x)$ for all $x\in W$;
\item\label{lm:adopt_quasi_tube:enum:crs_is_constant}
$\folEmb(c,x) = \crs(x)$ for all $x\in W$, i.e. $\folEmb^{-1}(\crs(W)) = c \times W$;
\item\label{lm:adopt_quasi_tube:enum:fixed_on_0}
if $\crs(0) = (c,0)$, then $\folEmb(t,0) = (t,0)$ for all $t\in(a,b)$.
\end{enumerate}
\end{lemma_eng}
\begin{proof}
Let $\pi:(a,b)\times \bR^n \to \bR^n$ be the standard projection.
Then the assumption that $\crs$ is a cross section means that the composition
\[
\pi \circ \crs: W \xrightarrow{~~\crs~~} (a,b)\times \bR^n \xrightarrow{~~\pi~~} \bR^n
\]
is an injective map between open subsets of $\bR^n$.
Hence by Brouwer's theorem on domain invariance, e.g.~\cite{HurewiczWallman:DimTh:1941}, $\pi\circ\crs(W)$ is an open neighborhood of $0$ in $\bR^n$.
Therefore we get an open embedding
\begin{align*}
\folEmb&: (a,b)\times W \to (a,b)\times \bR, &
\folEmb(t,x) = (t, \pi\circ\crs(x) )
\end{align*}
satisfying~\ref{lm:adopt_quasi_tube:enum:pres_lines} and~\ref{lm:adopt_quasi_tube:enum:fixed_on_0}.
Then $\folEmb^{-1}(\crs(W)) \subset (a,b)\times W$ can be regarded as a graph of certain continuous function $W \to (a,b)$.
Hence we get from Lemma~\ref{lm:make_constrant_graph} that $\folEmb$ can be composed with a foliated homeomorphism of $(a,b)\times W$ to satisfy~\ref{lm:adopt_quasi_tube:enum:crs_is_constant}, that is to make $\folEmb^{-1}(\crs(W))$ being the graph of the constant function $W \to c$.
Moreover, statements~\ref{lm:make_constrant_graph:enum:preserve_lines} and \ref{lm:make_constrant_graph:enum:fixed_on_0} of Lemma~\ref{lm:make_constrant_graph} allow to preserve properties~\ref{lm:adopt_quasi_tube:enum:pres_lines} and~\ref{lm:adopt_quasi_tube:enum:fixed_on_0} respectively.
\end{proof}

\begin{lemma_eng}\label{lm:extanding_foliated_chart}
Let $\leaf$ be a leaf of $\Partition$, $J_1, J_2 \subset \leaf$ be two compact segments such that $J_1 \cap J_2$ is a point, $\Vsp$ be an open $n$-disk, $\eps>0$, and
\begin{align*}
\folEmb_1& :(a-\eps,b+\eps)\times \Vsp \to \Xsp, &
\folEmb_2& :(b-\eps,c+\eps)\times \Vsp \to \Xsp
\end{align*}
be two open foliated embeddings such that
\begin{align*}
\folEmb_1\bigl( [a,b]\times 0\bigr) &= J_1, &
\folEmb_2\bigl( [b,c]\times 0\bigr) &= J_2, &
\folEmb_1(b,0) = \folEmb_2(b,0) = J_1 \cap J_2,
\end{align*}
and the union of the images of $\folEmb_1$ and $\folEmb_2$ does not contain compact leaves of $\Partition$.
Then there exists an open neighborhood $\Wsp$ of $0$ in $\Vsp$ and an open foliated embedding
\[
\folEmb:(a-\eps,c+\eps)\times \Wsp \to \Xsp
\]
such that $\folEmb\bigl([a,c]\times 0\bigr) = J_1\cup J_2$, see Figure~\ref{fig:join_fol_charts}.
\end{lemma_eng}
\begin{proof}
Notice that the assumption that the union of the images of $\folEmb_1$ and $\folEmb_2$ does not contain compact leaves of $\Partition$ implies that for any $\pa,\pb\in\Vsp$ the union of the arcs
\begin{align*}
&\folEmb_1\bigl( (a-\eps,b+\eps) \times \pa\bigr), &
&\folEmb_2\bigl( (b-\eps,c+\eps) \times \pb\bigr)
\end{align*}
does not contain a non-trivial loop, so the intersection of these arcs is connected (though possibly empty).

\begin{figure}[h]
\includegraphics[height=3cm]{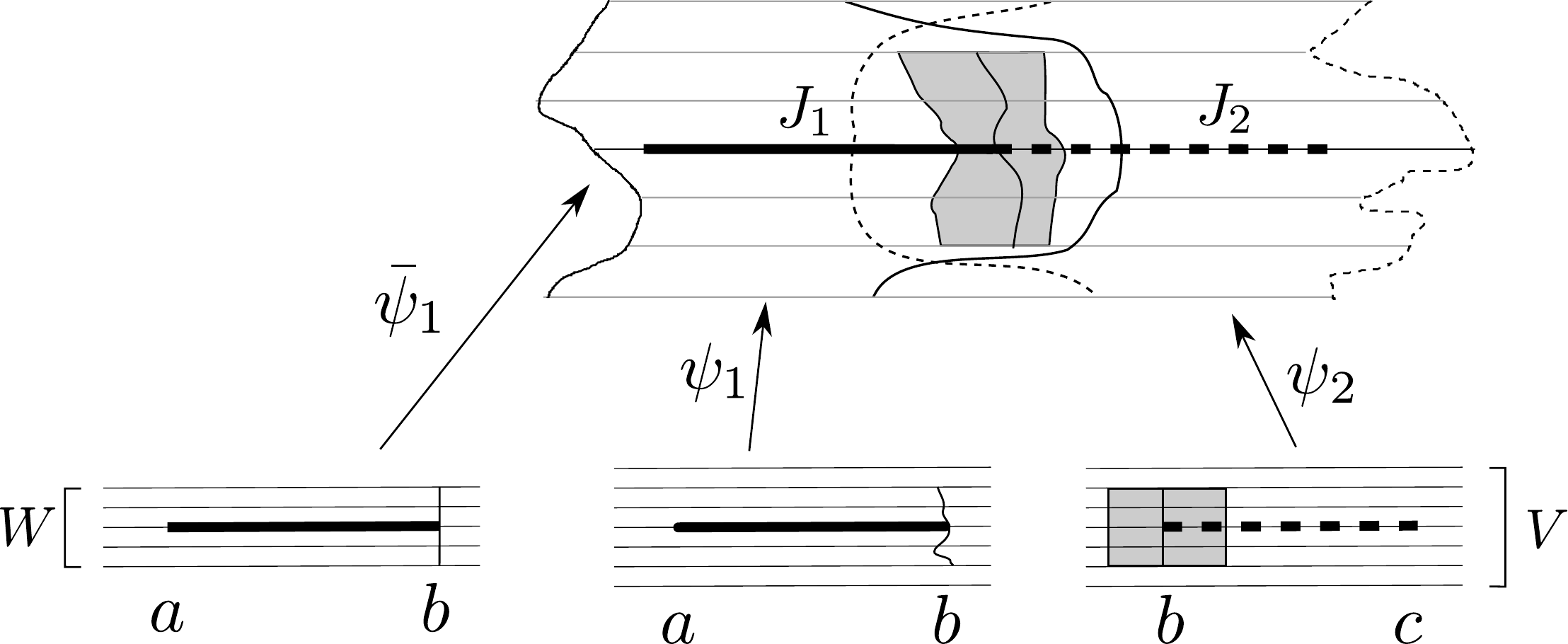}
\caption{}\label{fig:join_fol_charts}
\end{figure}

Since $\folEmb_1$ and $\folEmb_2$ are open embeddings, there exist $\delta>0$ and a small open neighborhood $\Wsp$ of $0$ in $\Vsp$ such that
\[ \folEmb_2\bigl( (b-\delta, b+\delta) \times W \bigr) \subset \mathrm{image}(\folEmb_1). \]
Then we have an embedding $\crs: W \to V$ defined by $\crs(\pa) = \folEmb_1^{-1}\bigl(\folEmb_2(b,\pa)\bigr)$, $\pa\in\Wsp$.
Hence by Lemma~\ref{lm:adopt_quasi_tube} one can find an open foliated embedding
\[
\bar{\folEmb}_1:(a-\eps,b+\eps)\times W \to \Xsp
\]
such that
\begin{itemize}
\item
$\bar{\folEmb}_1(\pt,0) = \folEmb_1(\pt,0)$ for all $\pt\in(a-\eps,b+\eps)$;
\item
$\mathrm{image}(\bar{\folEmb}_1) = \folEmb_1\bigl((a-\eps,b+\eps)\times W \bigr)$;
\item
$\pi\circ \bar{\folEmb}_1^{-1}(\px) = \pi \circ \folEmb_2^{-1}(\px)$ for all $\px\in \mathrm{image}(\bar{\folEmb}_1) \cap \mathrm{image}(\folEmb_2)$;
\item
$\bar{\folEmb}_1(b,\pa) = \folEmb_2(b,\pa)$ for all $\pa\in\Wsp$.
\end{itemize}
Now define the map $\folEmb :(a-\eps,c+\eps)\times \Vsp \to \Xsp$ by
\[
\folEmb(\pt,\pa) =
\begin{cases}
\bar{\folEmb_1}(\pt,\pa), & \pt\in(a-\eps, b], \\
\folEmb_2(\pt,\pa), & \pt\in[b, c+\eps).
\end{cases}
\]
One easily checks that $\folEmb$ is an open foliated embedding which coincides with $\folEmb_1$ on $(a-\eps,b]\times 0$ and with $\folEmb_2$ on $[b,c+\eps)\times0$.
In particular, $\folEmb\bigl([a,c]\times 0\bigr) = J_1\cup J_2$.
\end{proof}

\begin{corollary_eng}\label{cor:extanding_foliated_chart}{\rm c.f.~\cite[Lemma~22]{Meigniez:TrAMS:2002}}
Let $\Btrans$ be an open $n$-disk, $\leaf$ be a leaf of $\Partition$, and $J\subset \leaf$ be a compact segment.
Then there exists an open foliated embedding $\folEmb: (0,3) \times \Btrans \to \Xsp$ such that $\folEmb \bigl( [1,2] \times 0\bigr) = J$.
\end{corollary_eng}
\begin{proof}
Let us show that \myemph{there exists an open set $\Wsp$ such that $J \subset \Wsp$ and $\Wsp$ does not contain compact leaves of $\Partition$}.
Indeed, since $\leaf$ is either non-compact or is an embedded circle, it follows that $J \not= \leaf$.
Fix a point $x \in\leaf \setminus J$.
As $X$ is a regular space, there exist a pair of disjoint open neighborhoods $U_1 \ni x$ and $U_2 \supset J$.
Denote $W = U_2 \cap \SAT{U_1}$.
Then by Lemma~\ref{lm:saturation_is_open} $W$ is open.
Moreover, $J \subset \leaf = \SAT{x} \subset \SAT{U_1}$, so $J \subset W$.
Finally, since $W \subset \SAT{U_1}$, we see that $\SAT{y} \cap U_1 \neq \emptyset$ for each $y \in W$, that is $\SAT{y}\not\subset W$.
In other words, $W$ does not contain any leaf of $\Partition$.
In particular, $W$ can not contain compact leaves.

Notice that $J$ can be covered by finitely many foliated charts contained in $\Wsp$.
Lemma~\ref{lm:extanding_foliated_chart} allows to replace two consecutive foliated charts with one.
Hence the proof follows from that lemma by induction on the number of foliated charts covering $J$.
\end{proof}

\noindent
{\bf Cross sections.}
The following two lemmas describe general properties of cross sections.
\begin{lemma_eng}\label{lm:good_foliated_chart}
Let $\folEmb: (a,b) \times \Btrans \to \Xsp$ be an open foliated embedding.
Let also $\Usp = \folEmb\bigl( (a,b)\times \Btrans \bigr)$, $P_{\pa} = \folEmb\bigl( (a,b)\times \pa \bigr)$, $\pa\in\Btrans$, be a plaque of $\folEmb$, and $\crs:\crsDomain\to\Xsp$ be a cross section.
Then the following statements hold true.
\begin{enumerate}[leftmargin=*, label=$\rm(\arabic*)$]
\item\label{item:full_plaques}
Suppose $P_{\pa} \cap \crs(\crsDomain) \not=\varnothing$ for each $\pa\in\Btrans$.
Then then for each $\ps\in(a,b)$ the restriction map
\begin{equation}\label{eq:transverse_is_cross_section}
\folEmb|_{ \{\ps\} \times \Btrans }:  \{\ps\} \times \Btrans \to \Xsp
\end{equation}
is a cross section of $\Partition$.

\item\label{item:transverse_is_cross_section}
Suppose $\crs(\pb) \in P_{\pa}$ for some $\pa\in\Btrans$ and $\pb\in\crsDomain$.
Then there exists an open neighborhood $\crsDomain_{\pb}$ of $\pb$ in $\crsDomain$ and an open neighborhood $\Wsp_{\pb}$ of $\pb$ in $\Btrans$ such that
\begin{itemize}
\item
$\crs(\crsDomain_{\pb}) \subset \folEmb\bigl( (a,b) \times \Wsp_{\pb} \bigr)$;
\item
$P_{\pc} \cap \crs(\crsDomain_{\pb}) \not=\varnothing$ for each $\pc\in\Wsp_{\pb}$.
\end{itemize}
In particular, the restriction $\folEmb|_{ \ps\times \Wsp_{\pb}}:\ps\times \Wsp_{\pb}\to \Xsp$ is a cross sections of $\Partition$.

\item\label{item:cross_sections_at_each_point}
For every $\px \in \SAT{\crs(\crsDomain)}$ there exist an open subset $\Wsp$ of $\bR^n$ and a cross-section $\folEmb_{\px}: \Wsp \to \Xsp$ such that $\px \in \folEmb_{\px}(\Wsp) \subset \SAT{\crs(\crsDomain)}$.
\end{enumerate}
\end{lemma_eng}
\begin{proof}
\ref{item:full_plaques}
Suppose $P_{\pa} \cap \crs(\crsDomain) \not=\varnothing$ for all $\pa\in\Btrans$.
Since $\crs(\crsDomain)$ intersects each leaf of $\Partition$ in at most one point, it follows that distinct plaques $P_{\pa}$ and $P_{\pb}$ for $\pa\not=\pb\in\Btrans$ belong to distinct leaves of $\Partition$.
As $\Btrans$ is an open subset of $\bR^n$, the map~\eqref{eq:transverse_is_cross_section} is a cross section for each $\ps\in (a,b)$.

\medskip

\ref{item:transverse_is_cross_section}
Consider the following map:
\begin{equation}\label{eq:open_map}
\xi = \pi \circ \folEmb \circ \crs|_{\crs^{-1}(\Usp)}:
 \crs^{-1}(\Usp) \ \xrightarrow{~~\crs~~} \
 \Usp  \ \xrightarrow{~~\folEmb~~} \
 (a,b)\times \Btrans  \ \xrightarrow{~~\pi~~} \  \Btrans,
\end{equation}
where $\pi$ is the standard projection to the second coordinate.

Then the assumption $\crs(\pb) \in P_{\pa}$ for some $\pa\in\Btrans$ and $\pb\in\crsDomain$ implies that $\pb\in\crs^{-1}(\Usp)$ and $\xi(\pb) = \pa$.

Since the images of distinct points of $\crsDomain$ under $\crs$ are contained in distinct leaves of $\Partition$, they also belong to distinct plaques of $\folEmb$, whence $\xi$ is an injective continuous map between open subsets of $\bR^n$.
Hence, by Brouwer theorem on domain invariance $\xi$ is an open map, \cite{Brower:1912}.
In particular, $\xi$ yields a homeomorphism of some open neighborhood $\crsDomain_{\pb}$ of $\pb$ onto some open neighborhood $\Wsp_{\pa}$ of $\pa$ in $\Btrans$.
This implies that $\crs(\crsDomain_{\pb}) \subset \folEmb\bigl( (a,b) \times \Wsp_{\pa} \bigr)$ and $P_{\pc} \cap \crs(\crsDomain_{\pb}) \not=\varnothing$ for each $\pc\in\Wsp_{\pa}$.

\medskip

\ref{item:cross_sections_at_each_point}
Let $\leaf$ be the leaf containing $\px$ and $\py = \crs(\pb) = \leaf \cap \crs(\crsDomain)$.
If $\px=\py$, then one can put $\Wsp_{\px}=\crsDomain$ and $\crs_{\px} = \crs$.

Therefore suppose $\px\not=\py$.
Let $J \subset \omega$ be a closed segment with ends $\px$ and $\py$.
Then by Corollary~\ref{cor:extanding_foliated_chart} there exists an open foliated embedding $\folEmb: (0,3) \times \Btrans \to \Xsp$ such that $\folEmb \bigl( [1,2] \times 0\bigr) = J$, $\folEmb(1,0)=\px$ and $\folEmb(2,0) = \py$.

Thus $\py = \crs(\pb) = \folEmb(2,0) \in P_{0} = \folEmb\bigl( (0,3) \times 0 \bigr)$, and so by~\ref{item:transverse_is_cross_section} there exists a neighborhood $\Wsp$ of $0$ in $\Btrans$ such that the map $\folEmb_x:\Wsp\to\Xsp$ defined by $\folEmb_x(\pc) = \folEmb(1, \pc)$ is a cross section with $\folEmb_{\px}(\Wsp) \subset \SAT{\crs(\crsDomain)}$.
It remains to note that $\folEmb_x(0)=\folEmb(1,0)=\px$.
\end{proof}

\begin{lemma_eng}\label{lemma:proper_leaf1}
Let $\leaf_0$ be a leaf of $\Partition$.
Suppose that for each leaf $\leaf$ of $\Partition$ contained in $\SAT{\Cl{\leaf_0}}$ there exists a cross section $\crs$ passing through $\leaf$.
Then $\leaf_0$ is properly embedded, i.e. it satisfyies conditions (m) and (c) of Lemma~\ref{lm:proper_leaves}.
\end{lemma_eng}
\begin{proof}
If $\leaf_0$ is compact, then it is necessarily properly embedded.
Therefore assume that $\leaf_0$ is non-compact.

(m)
Let $\leaf \subset \SAT{\Cl{\leaf_0}}$ be a leaf of $\Partition$.
By \ref{item:cross_sections_at_each_point} of Lemma~\ref{lm:good_foliated_chart} for each $\px \in \leaf$ there exists an open foliated embedding $\folEmb:(-1,1)\times \Btrans \to \Xsp$ such that $\folEmb(0,0)=\px$ and different plaques of $\folEmb$ are contained in different leaves of $\Partition$.
In particular, $\folEmb$ homeomorphically maps $(-1,1)\times\{0\}$ onto an open neighbourhood of $\px$ in $\leaf$.
This implies that $\leaf$ is an embedded 1-submanifold of $\Xsp$.

(c)
Let $\px \in \SAT{\Cl{\leaf_0}} \setminus \leaf_0$.
Then decreasing $\Btrans$ one can assume that the image of $\folEmb$ does not intersect $\leaf_0$, whence $x \notin \Cl{\leaf_0}$.
From arbitrariness of $x \in \SAT{\Cl{\leaf_0}}$ we conclude that $\leaf_0$ is closed in $\Xsp$.
\end{proof}

\noindent
{\bf Parallel cross sections.}
Let $\crs:\crsDomain\to\Xsp$ be a cross section and $\crsSubDomain \subset \crsDomain$ be an open subset.
Then a cross section $\delta:\crsSubDomain \to \Xsp$ \myemph{parametrically agrees} with $\crs$, whenever for each $\pa\in\crsSubDomain$ the points $\delta(\pa)$ and $\crs(\pa)$ belong to the same leaf.
Also $\delta$ is \myemph{parallel to $\crs$} if it parametrically agrees with $\crs$ and $\delta(\crsSubDomain) \cap \crs(\crsSubDomain) = \varnothing$.

Let $\crs_0,\crs_1:\crsDomain \to \Xsp$ be two parallel cross sections intersecting only non-compact leaves.
For each $\pa\in\crsDomain$ let $\pleaf{\pa}$ be the leaf containing $\gamma_0(\pa)$ and $\gamma_1(\pa)$, $\segm{\pa} \subset \pleaf{\pa}$ be the compact segment with ends $\gamma_0(\pa)$ and $\gamma_1(\pa)$, and $\Int{\segm{\pa}}$ be the interior of $\segm{\pa}$.
In this situation we will put:
\begin{align}
\Lset{\gamma_0}{\gamma_1} &:= \mathop{\cup}\limits_{\pa\in\crsDomain} \Int{\segm{\pa}}, &
\Kset{\gamma_0}{\gamma_1} &:= \mathop{\cup}\limits_{\pa\in\crsDomain} \segm{\pa}.
\end{align}

\begin{figure}[h]
\includegraphics[height=1.2cm]{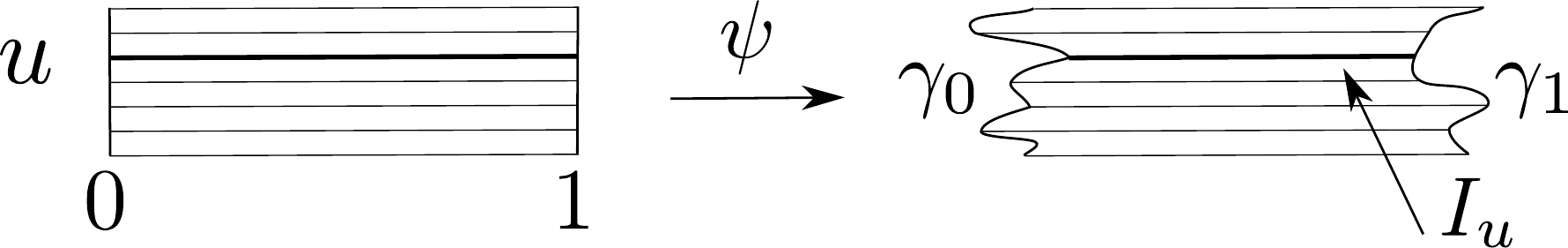}
\caption{}\label{fig:space_between_crs}
\end{figure}

\begin{lemma_eng}\label{lm:parallel_crossections}
There exists a homeomorphism $\folEmb:[0,1]\times\crsDomain \to \Kset{\gamma_0}{\gamma_1}$ such that
\begin{align*}
\folEmb\bigl([0,1]\times \pa\bigr) &= \segm{\pa}, &
\folEmb(0,\pa) &= \crs_0(\pa), &
\folEmb(1,\pa) &= \crs_1(\pa)
\end{align*}
for every $\pa\in\crsDomain$, see Figure~\ref{fig:space_between_crs}.
In particular, $\folEmb\bigl( (0,1)\times\crsDomain \bigr) = \Lset{\gamma_0}{\gamma_1}$.
Moreover, $\Lset{\gamma_0}{\gamma_1}$ is open in $\Xsp$.
\end{lemma_eng}
\begin{proof}
Fix some $\eps>0$ and denote $J = (-\eps, 1+\eps)$.
Then it follows from Corollary~\ref{cor:extanding_foliated_chart} and Lemma~\ref{lm:adopt_quasi_tube} that for each $\pa\in\crsDomain$ there exists a neighborhood $\Wsp_{\pa}$ in $\crsDomain$ and an open foliated embedding $\folEmb_{\pa}: J \times{\Wsp_{\pa}}\to \Xsp$ having the following properties:
\begin{enumerate}[leftmargin=2em]
\item[(a)]
$\folEmb_{\pa}([0,1]\times\pa) = \segm{\pa}$, \ $\folEmb_{\pa}(0,\pa) = \gamma_0(\pa)$, \ and \ $\folEmb_{\pa}(1,\pa) = \gamma_1(\pa)$;
\item[(b)]
$\folEmb_{\pa}(J\times\pb)$, $\gamma_0(\pb)$, and therefore $\gamma_1(\pb)$, are contained in the same leaf of $\Partition$ for each $\pb\in{\Wsp_{\pa}}$;
\item[(c)]
$\gamma_0({\Wsp_{\pa}}) = 0 \times {\Wsp_{\pa}}$ \ and \ $\gamma_1({\Wsp_{\pa}}) = 1 \times {\Wsp_{\pa}}$.
\end{enumerate}
In particular, this implies that the set
\[
\Lset{\gamma_0}{\gamma_1} = \mathop{\cup}\limits_{x\in\crsDomain} \, \folEmb_{\pa}\bigl( (0,1)\times \Wsp_{\pa} \bigr)
\]
is open in $\Xsp$.

As $\crsDomain$ is paracompact, there is a locally finite cover $\{\Wsp_{i}\}_{i\in\Lambda}$ of $\crsDomain$ and for each $i\in\Lambda$ an open foliated embedding $\folEmb_i:J\times \Wsp_{i} \to \Xsp$ such that $\folEmb_i([0,1]\times\pa) = \segm{\pa}$ for all $\pa\in \Wsp_{i}$.
Denote $\Usp_{i} = \folEmb_i(J\times \Wsp_{i})$ and $\Usp = \mathop{\cup}\limits_{i\in\Lambda} \Usp_{i}$.
Then $\Usp$ is an open neighborhood of $\Kset{\gamma_0}{\gamma_1}$ and $\{\Usp_{i}\}_{i\in\Lambda}$ is a locally finite cover of $\Usp$.

Let $\{\lambda_i:\crsDomain\to[0,1]\}_{i\in\Lambda}$ be a partition of unity subordinated to the cover $\{\Wsp_{i}\}_{i\in\Lambda}$.
Thus $\mathrm{supp}(\lambda_i) \subset \Wsp_{i}$ and $\sum_{i\in\Lambda} \lambda_i(\pa) = 1$.
Let also $p_i:J\times \Wsp_{i} \to J$ and $q_i:J\times \Wsp_{i} \to \Wsp_{i}$ be the standard projections, and
\[
 \mu_i = \lambda_i \circ q_i \circ \folEmb_i^{-1}:
 \ \Usp_{i} \ \xrightarrow{~~\folEmb_i^{-1}~~} \ J\times \Wsp_{i} \ \xrightarrow{~~q_i~~} \ \Wsp_{i} \ \xrightarrow{~~\lambda_i~~} \ [0,1].
\]
Then $\mathrm{supp}(\mu_i) = J\times\mathrm{supp}(\lambda_i)$, whence $\mu_i$ extends by zero to a continuous function on all of $\Usp$.

Let $f: \Usp \to J$ be the function defined by the following rule:
\[f(\px) = \sum_{\px\in \Usp_{i}} \mu_{i}(\px) \ \cdot \ p_{i}\circ \folEmb_i^{-1}(\px).\]
Since for each $\pa\in \Wsp_{i}$ the function $p_{i}\circ \folEmb_i^{-1}: \segm{\pa} \to [0,1]$ is homeomorphism which maps $\gamma_0(\pa)$ and $\gamma_1(\pa)$ to $0$ and $1$ respectively, and $\sum_{j\in\Lambda} \mu_j\equiv 1$, we see that the restriction $f|_{\segm{\pa}}$ is a convex linear combination of orientation preserving homeomorphisms.
Therefore $f|_{\segm{\pa}}: \segm{\pa}\to[0,1]$ is a homeomorphism as well.

Let also $g:\Usp \to \crsDomain$ be the map defined by $g(\px) = q_i(\px)$ whenever $\px\in \Usp_{i}$.
Due to (b) this definition does not depend on a particular $\Usp_{i}$ containing $\px$.
Hence $g$ is a well-defined continuous map.

Then the mapping \[ \phi = (f, g): \Kset{\gamma_0}{\gamma_1} \to [0,1]\times\crsDomain\] is a continuous bijection being also a local homeomorphism, and so it is a homeomorphism.
Moreover, $\phi(\segm{\pa}) = [0,1]\times \pa$ for all $\pa\in\crsDomain$.
Therefore $\folEmb = \phi^{-1}$ is the required homeomorphism.
\end{proof}

\begin{lemma_eng}\label{lm:tube_along_leaf}
Let $\crs_i:\crsDomain\to\Xsp$, $i\in\bZ$, be a family of pairwise parallel cross sections intersecting only non-compact leaves and $\Usp = \SAT{\crs_i(\crsDomain))}$ be the common saturation of their images.
Suppose also that the following two conditions hold:
\begin{enumerate}[leftmargin=2em]
\item[\rm(1)]
$\Lset{\crs_{i}}{\crs_{i+1}} \cap \Lset{\crs_{j}}{\crs_{j+1}} = \varnothing$ for $i\not=j$;
\item[\rm(2)]
$\mathop{\cup}\limits_{i=-\infty}^{\infty} \Kset{\crs_{i}}{\crs_{i+1}} = \Usp$.
\end{enumerate}
Then $\Usp$ is open in $\Xsp$ and foliated homeomorphic with $\bR\times\crsDomain$.
\end{lemma_eng}
\begin{proof}
By Lemma~\ref{lm:parallel_crossections} for each $i\in\bZ$ there exists a homeomorphism
\[\folEmb_{i}:[i,i+1] \times \crsDomain \to \Kset{\crs_{i}}{\crs_{i+1}}\]
such that for each $\pa\in\crsDomain$
\begin{enumerate}[leftmargin=2em, label=$\bullet$]
\item
$\folEmb_{i}\bigl([i,i+1]\times \pa)$ is a segment of the leaf of $\Partition$ between the points $\crs_i(\pa)$ and $\crs_{i+1}(\pa)$;
\item
$\folEmb_{i-1}(i,\pa) = \folEmb_{i}(i,\pa) = \crs_i(\pa)$.
\end{enumerate}
Therefore we have a homeomorphism
\[ \folEmb: \bR \times \crsDomain \longrightarrow \mathop{\cup}\limits_{i=-\infty}^{\infty} \Kset{\crs_{i}}{\crs_{i+1}} = \Usp \]
defined by $\folEmb(\pt,\pa) = \folEmb_{i}(\pt,\pa)$ whenever $\pt\in[i,i+1]$ and $\pa\in\crsDomain$.
Moreover, $\Usp = \mathop{\cup}\limits_{i=1}^{\infty} \Lset{\crs_{-i}}{\crs_{i}}$ is open in $\Xsp$.
\end{proof}

\section{Proof of Theorem~\ref{th:trivialization_of_cross_section}}\label{sect:proof_of_theorem:th:trivialization_of_cross_section}
Let $\crs:\crsDomain \to \Xsp$ be a cross section intersecting only leaves being simultaneously non-compact, properly embedded, and non-special.
We have to prove that its saturation $\SAT{\crs(\crsDomain)}$ is open and foliated homeomorphic with $\bR\times\crsDomain$.

First we will assume that $\partial\Xsp = \varnothing$.
The proof of the case $\partial\Xsp\not=\varnothing$ will follow from the case $\partial\Xsp=\varnothing$ by passing to the double $2\Xsp$ of $\Xsp$ and considering the one-dimensional foliation on $2\Xsp$ induced by $\Partition$.
It will be given at the end op this section.

Our proof is based on the following statement which will be proved below.
\begin{proposition_eng}\label{pr:global_crs_outside_compact_set}
Let $K \subset \Xsp$ be a compact subset.
Then one can find two parallel cross sections $\alpha,\beta:\crsDomain\to\Xsp$ parametrically agreeing with $\crs$ and satisfying
\[
\SAT{\crs(\crsDomain)} \cap K \ \subset \ \Lset{\alpha}{\beta}.
\]
Moreover, if $\acrs,\bcrs:\crsDomain\to\Xsp$ are two parallel cross sections parametrically agreeing with $\crs$, then one can assume that
\[
\bigl( \SAT{\crs(\crsDomain)} \cap K \bigr) \ \bigcup \ \Kset{\acrs}{\bcrs} \ \subset \ \Lset{\alpha}{\beta}.
\]
\end{proposition_eng}

Before proving Theorem~\ref{th:trivialization_of_cross_section} let us deduce it from Proposition~\ref{pr:global_crs_outside_compact_set}.

\begin{figure}[h]
\includegraphics[height=2cm]{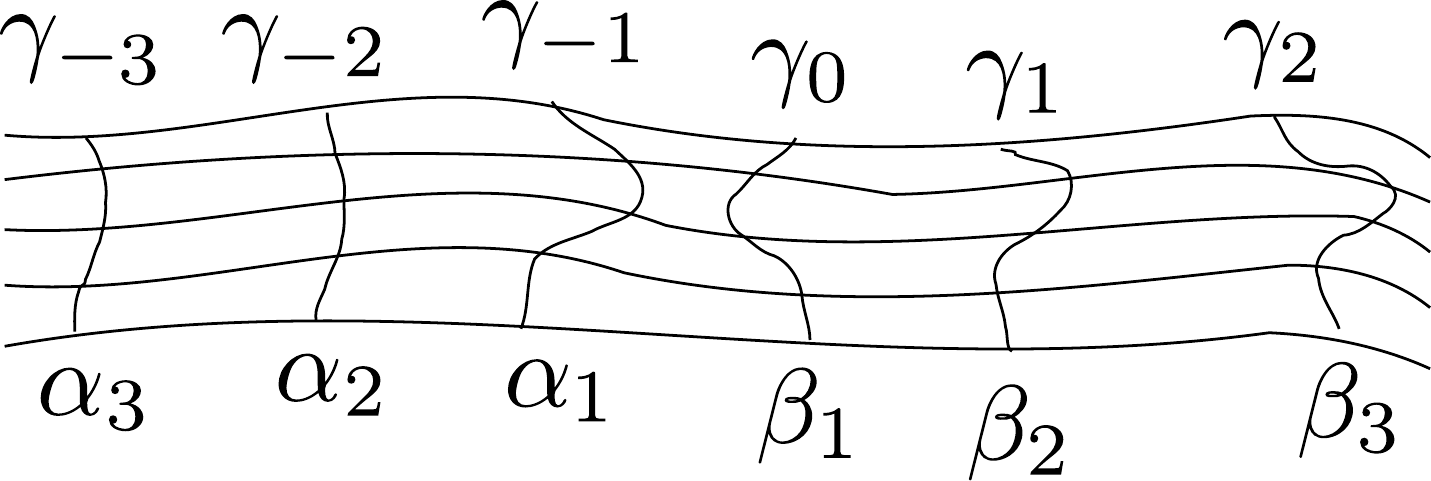}
\caption{}\label{fig:inf_cross_sect_sequence}
\end{figure}

Fix any increasing sequence $K_1\subset K_2\subset\cdots$ of compact subsets of $\Xsp$ such that $\Xsp = \mathop{\cup}\limits_{i\in\bN} K_i$.
Using Proposition~\ref{pr:global_crs_outside_compact_set} one constructs a family of parallel cross sections $\alpha_i,\beta_i:\crsDomain \to \Xsp$, $i\in\bN$, parametrically agreeing with $\crs$ and such that
\begin{enumerate}[leftmargin=2em]
\item[\rm(1)]
$\SAT{\crs(\crsDomain)} \cap K_i \ \subset \ \Lset{\alpha_i}{\beta_i}$;
\item[\rm(2)]
$\Kset{\alpha_{i-1}}{\beta_{i-1}} \  \subset \ \Lset{\alpha_{i}}{\beta_{i}}$ for all $i\geq2$.
\end{enumerate}
Hence
\[
\SAT{\crs(\crsDomain)} \ = \ \mathop{\cup}\limits_{i\in\bN} \SAT{\crs(\crsDomain)}\,\cap\,K_i \ = \ \mathop{\cup}\limits_{i\in\bN} \Lset{\alpha_i}{\beta_i}.
\]
Exchanging $\alpha_i$ and $\beta_i$ if necessary and re-denoting them as follows: $\crs_{-i} = \alpha_i$, and $\crs_{i-1} = \beta_i$ for $i\in\bN$, one can assume that the sequence of cross sections $\{\crs_i\}_{i\in\bZ}$ satisfies assumptions of Lemma~\ref{lm:tube_along_leaf}, see Figure~\ref{fig:inf_cross_sect_sequence}.
Hence $\SAT{\crs(\crsDomain)}$ is open and foliated homeomorphic with $\bR\times\crsDomain$.
This proves Theorem~\ref{th:trivialization_of_cross_section} modulo Proposition~\ref{pr:global_crs_outside_compact_set}.

\medskip

The following lemma guarantees existence of local cross sections in Proposition~\ref{pr:global_crs_outside_compact_set}.

\begin{lemma_eng}\label{lm:local_crs_outside_compact_set}
Let $K \subset \Xsp$ be a compact subset.
Then for each $\pa\in\crsDomain$ one can find an open neighborhood $\crsSubDomain$ in $\crsDomain$ and two parallel cross sections $\alpha,\beta:\crsSubDomain\to\Xsp$ parametrically agreeing with $\crs$ and such that
\[  \SAT{\crs(\crsSubDomain)} \cap K \ \subset \ \Lset{\alpha}{\beta}. \]
\end{lemma_eng}
\begin{proof}
Suppose that lemma fails, so there exists $\pa\in\crsDomain$ belonging to some leaf $\leaf$ such that
\begin{itemize}
\item
for any decreasing sequence $W_i$ of neighborhoods of $\pa$ in $\crsDomain$ with $\mathop{\cap}\limits_{i\in\bN} \crsSubDomain_i = \{\pa\}$
\item
and any family of pairs of parallel cross sections $\alpha_i,\beta_i: \crsSubDomain_i \to \Xsp$, $i\in\bN$, parametrically agreeing with $\crs$
\end{itemize}
the set
\[
\SAT{\crs(\crsSubDomain_i)}\setminus \Lset{\alpha_i}{\beta_i}
\]
contains some point $\px_{i}\in K$.

Denote
\[
\Usp = \mathop{\cup}\limits_{i\in\bN} \Lset{\alpha_i}{\beta_i}.
\]
Then one can assume, in addition, that the following properties hold:
\begin{enumerate}[leftmargin=2em]
\item[(a)]
the sequence $\{\px_{i}\}_{i\in\bN}$ converges to some point $\px\in K$;
\item[(b)]
$\leaf \subset \Usp$;
\item[(c)]
$\px_{i}\not\in\Usp$ for all $i\in\bN$, whence $\px\not\in\Usp$ as well, and so $\px\not\in\leaf$.
\end{enumerate}
Indeed, (a) follows from compactness of $K$.

To prove (b) fix any continuous bijection $\bijleaf:\bR\to\leaf$.
By assumption $\bijleaf$ is proper, so one can find $A>0$ such that $\leaf\cap K \subset \bijleaf(-A,A)$.
Choose $\alpha_i$ and $\beta_i$ so that $\alpha_i(\Wsp_i) \cap K = \beta_i(\Wsp_i) \cap K = \varnothing$,
\begin{multline*}
\cdots < \alpha_{i+1}(\pa)  < \alpha_{i}(\pa) < \cdots < \alpha_{1}(\pa) < -A < \\ < A < \beta_{1}(\pa) < \cdots < \beta_{i}(\pa) < \beta_{i+1}(\pa)  < \cdots
\end{multline*}
$\lim\limits_{i\to+\infty} \alpha_{i}(\pa)  = -\infty$, and $\lim\limits_{i\to+\infty}\beta_{i}(\pa)  = +\infty$.
Then we will have that $\leaf \subset \Usp$.

Finally, to satisfy (c) choose $\crsSubDomain_{i+1}$ so small that $\px_{i}\not\in \SAT{\crs(\crsSubDomain_{i+1})}$ for all $i\in\bN$.

Now let $\leaf_i$ be the leaf of $\Partition$ containing $\px_{i}$, and $\py_{i} = \leaf_i \cap \crs(\crsSubDomain_1)$.
Then the sequence $\{\py_{i}\}_{i\in\bN}$ converges to $\py =\crs(\pa) \in \leaf$.
Hence $\prj(\px) \not = \prj(\py) = \prj(\leaf)$, while $\prj(\px_{i}) = \prj(\py_{i}) = \prj(\leaf_i)$ for all $i\in\bN$.
Therefore by Lemma~\ref{lm:charact_spec_leaves} $\leaf$ is a special leaf which contradicts to the assumption.
\end{proof}

The rest of the proof of Theorem~\ref{th:trivialization_of_cross_section} is based on E.~Michael's result about selections, \cite{Michael:AM:2:1956}.

Let $2^{\Xsp}$ be the set of all subsets of $\Xsp$ and $\mathcal{E}(\Xsp) \subset 2^{\Xsp}$ be the set of all closed subsets of $\Xsp$.
Let also $A\subset\crsDomain$ be a subset and $q:\crsDomain \Rightarrow \Xsp$ be a multivalued map, i.e. a map $q:\crsDomain \to 2^{\Xsp}$.
Then a \myemph{selection} for the restriction $q|_{A}$ is a continuous map $\phi:A \to \Xsp$ such that $\phi(\px) \in q(\px)$ for all $\px\in A$.

A multivalued map $q:\crsDomain\Rightarrow \Xsp$ is called \myemph{lower semi-continuous} whenever for each open $\Usp \subset \Xsp$ the set
\[ T_{\Usp} := \{\px\in\crsDomain \mid q(\px) \cap \Usp \not=\varnothing\} \]
 is open in $\crsDomain$.

A family $\Family \subset 2^{\Xsp}$ is called \myemph{equi-$LC^{k}$}, $k\geq0$, if for every $\Pel\in\Family$, $\px\in\Pel$, and a neighborhood $\Usp_{\px}$ of $\px$ in $\Xsp$, there exists a neighborhood $O_{\px}$ of $\px$ in $\Xsp$ such that for every $\Qel \in \Family$ every continuous map $f:S^m \to \Qel\cap O_{\px}$ of an $m$-sphere ($m\leq k$) is homotopic to a constant map in $\Qel\cap U_{\px}$.

A topological space $\Zsp$ is called $C^k$, or \myemph{$k$-connected}, $k\geq0$, if every continuous map $f:S^m \to \Zsp$ of an $m$-sphere ($m\leq k$) is homotopic to a constant map.

\begin{theorem_eng}\label{th:MichaelTheorm}{\rm\cite[Theorem~1.2]{Michael:AM:2:1956}}
Let $\crsDomain$ be a separable metric space, $A \subset \crsDomain$ be a closed subset with $\dim (\crsDomain\setminus A) \leq k+1$, $\Xsp$ a complete metric space, $\Family \subset \mathcal{E}(\Xsp)$ be equi-$LC^k$ and $q:\crsDomain \to \Family$ be a lower semi-continuous map.
Then every selection for $q|_{A}$ can be extended to a selection for $q|_{U}$ for some open $U \supset A$.
If also every $S \subset \Family$ is $C^k$, then one can take $U = X$.
\end{theorem_eng}

We will use the following particular case of Theorem~\ref{th:MichaelTheorm}.
\begin{corollary_eng}\label{cor:MichaelTheorm}
Let $\crsDomain$ be a separable metric space, $\dim\crsDomain = n$, $\Xsp$ be a complete metric space, and $\Family \subset \mathcal{E}(\Xsp)$ be an equi-$LC^{n+1}$ family such that each $\Qel\in\Family$ is contractible.
Then every lower semi-continuous multivalued map $q:\crsDomain \to \Family$ has a continuous selection.
\end{corollary_eng}

\medskip 
\noindent
{\bf Proof of Proposition~\ref{pr:global_crs_outside_compact_set}}
Since $\crsDomain$ is paracompact, it follows from Lemma~\ref{lm:local_crs_outside_compact_set} that there exist
\begin{itemize}
\item
a locally finite open cover $\mathcal{W} = \{\crsSubDomain_i\}_{i\in\bN}$ of $\crsDomain$ with compact closures $\overline{\Wsp_i}$, and
\item
a family of pairs of parallel cross sections $\alpha_i,\beta_i:\crsSubDomain_i\to\Xsp$, $i\in\bN$, parametrically agreeing with $\crs$
\end{itemize}
such that
\[
\SAT{ \crs(\crsSubDomain_i) } \cap  K_i \ \subset \ \Lset{\alpha_i}{\beta_i},
\]
where $K_i = K \cup \acrs(\overline{\Wsp_i}) \cup \bcrs(\overline{\Wsp_i})$ whenever the cross sections $\acrs,\bcrs:\crsDomain\to\Xsp$ are given, and $K_i = K$ otherwise.

Then it follows from Lemmas~\ref{lm:parallel_crossections} and~\ref{lm:make_constrant_graph} that for each $i\in\bN$ one can find an embedding $\folEmb_{i}: [-1,1] \times \crsSubDomain_i \to \Xsp$ such that for each $\pa\in\crsDomain$
\begin{enumerate}[leftmargin=2em]
\item[(1)]
$\folEmb_{i}\bigl( [-1,1] \times \pa \bigr)$ is contained in the leaf of $\Partition$;
\item[(2)]
$\folEmb_{i}(-1,\pa) = \alpha_i(\pa)$, \ $\folEmb_{i}(0,\pa) = \crs(\pa)$, \ $\folEmb_{i}(1,\pa) = \beta_i(\pa)$;
\item[(3)]
$\SAT{\crs(\crsSubDomain_i)} \ \cap \  K_i \ \subset \ \folEmb_{i}\bigl( (-1,1) \times \crsSubDomain_i \bigr)$;
\item[(4)]
$\alpha_i(\crsSubDomain_i)$ are contained in the same path component of $\SAT{\crs(\crsDomain)} \setminus \crs(\crsDomain)$ for all $i\in\bN$.
\end{enumerate}

Let $\pa\in\crsDomain$, $\pleaf{\pa}$ be the leaf of $\Partition$ containing $\crs(\pa)$, and $\phi_{\pa}:\bR\to\pleaf{\pa}$ be any bijection satisfying $\phi_{\pa}^{-1}(\alpha_i(\pa))<0$, and $\phi_x(0)=\crs(\pa)$.
Therefore $\phi_{\pa}^{-1}(\beta_i(x))>0$ for all $i$ such that $\pa\in\crsSubDomain_i$.
Then there are two numbers $a_{\pa}, b_{\pa}$ such that
\[
\pleaf{\pa} \ \setminus \ \mathop{\cap}\limits_{i \ : \ \pa\in \crsSubDomain_i} \Lset{\alpha_i}{\beta_i}.
\]
consists of two half closed intervals $A_{\pa} = \phi_{\pa}(-\infty,a_{\pa}]$ and $B_{\pa} = \phi_{\pa}[b_{\pa},+\infty)$.

Since $\pleaf{\pa}$ is a properly embedded leaf, it follows that $A_{\pa}$ and $B_{\pa}$ are closed in $\Xsp$.
Moreover, by (3) they do not intersect $K$.

Define the following two maps $a, b: \crsDomain \to \mathcal{E}(\Xsp)$, i.e.\! multivalued mappings $a,b:\crsDomain\Rightarrow\Xsp$ with closed images, by
\begin{align*}
a(\pa) &= A_{\pa}, &
b(\pa) &= B_{\pa}
\end{align*}
for $\pa\in\crsDomain$.

\begin{lemma_eng}\label{lm:ab_prop}
{\rm(i)}~The maps $a$ and $b$ are lower semi-continuous.

{\rm(ii)}~The families $\mathcal{A} = \{ A_{\pa} \mid \pa\in\crsDomain \}$ and $\mathcal{B} = \{ B_{\pa} \mid \pa\in\crsDomain \}$ are equi-$LC^k$ for all $k\geq0$.
\end{lemma_eng}
\begin{proof}
It suffices to check (i) and (ii) for $a$ only.

(i) We should check that for each open $U \subset \Xsp$ the set 
\[ T_{U} = \{\pa\in\crsDomain \mid a(\pa) \cap U \not=\varnothing\}\] is open as well.

Let $\pa\in\crsDomain$ be such that $A_{\pa} \cap U \not=\varnothing$, and $\px\in A_{\pa} \cap U$.
Since $U$ is open, one can assume that $\px$ is not the end of $A_{\pa}$, that is $\phi_{\pa}^{-1}(y) < a_{\pa}$.

By assumption $\pa \in W_i$ for some $i\in\bN$.
Then by Corollary~\ref{cor:extanding_foliated_chart} for the closed interval on $\pleaf{\pa}$ between $\px$ and $\phi_{\pa}(a_{\pa})$ there exists an open neighborhood $O$ of $\pa$ in $\Wsp_i$ and an open foliated embedding $\folEmb:(-1,2) \times O \to \Xsp$ such that
\begin{itemize}
\item[\rm(a)]
$\folEmb\bigl( (-1,2) \times O \bigr) \subset \Usp$;
\item[\rm(b)]
$\folEmb(0,\pa) = \px$;
\item[\rm(c)]
$\folEmb\bigl( (-1,2) \times \pb \bigr) \subset \pleaf{\pb}$;
\item[\rm(d)]
$\folEmb(1,\pb) = \alpha_i(\pb)$.
\end{itemize}
It follows from (a) and (d) that $\folEmb((-1,1] \times \pb) \subset A_{\pb} \cap \Usp$, whence $O \subset T_U$.
Thus $T_{U}$ is open, and so $a$ is a lower semi-continuous multivalued map.

\smallskip

(ii) Notice that for each $\px\in\Xsp$ there exists an open neighborhood $U_{\px}$ such that the intersection of $U_{\px}$ with each leaf $\leaf$ is either empty or homeomorphic to an open interval.
Therefore intersection of $U_{\px}$ with each set $A_{\pa}$ is either empty or homeomorphic to $(0,1)$ or to $(0,1]$.
In the latter two cases $U_{\px} \cap A_{\pa}$ is contractible.
Hence every continuous map $S^k\to U_{\px} \cap A_{\pa}$ is null homotopic and one can put $O_y = U_y$.
This means that $\mathcal{A}$ is equi-$LC^k$ for all $k\geq0$.
\end{proof}

Since for each $\pa\in\crsDomain$ the sets $A_{\pa}$ and $B_{\pa}$ are contractible, it follows from Lemma~\ref{lm:ab_prop} that $a$ and $b$ satisfy assumptions of Corollary~\ref{cor:MichaelTheorm}.
Hence they admit continuous selections $\alpha,\beta:\crsDomain\to\Xsp$ and these selections are the required cross sections.
This completes Proposition~\ref{pr:global_crs_outside_compact_set}.

\medskip

\noindent
{\bf Proof of Theorem~\ref{th:trivialization_of_cross_section}. Case $\partial\Xsp \not= \varnothing$.}
We need the following simple lemma whose proof we leave for the reader.
\begin{lemma_eng}\label{lm:double_of_open_set}
Let $\xi:\bR^n \to \bR^n$ be the involution defined by \[\xi(x_1,\ldots,x_{n-1},x_{n}) = (x_1,\ldots,x_{n-1},-x_{n}).\]
Then for each subset $\crsDomain \subset \bR^n_{+}$ open in the induced topology of $\bR^n_{+}$, its \myemph{double} $\dblV = \crsDomain\cup \xi(\crsDomain)$ is open in $\bR^n$.
\qed
\end{lemma_eng}

Now let
\[
\dblX = \Xsp_1 \mathop{\sqcup}\limits_{\mathrm{id}:\partial\Xsp_1\to\partial\Xsp_2} \Xsp_2
\]
be the double of $\Xsp$, i.e.\! the union of two copies $\Xsp_1$ and $\Xsp_2$ of $\Xsp$ glued along their boundaries by the identity map.
Let also $\sigma:\dblX\to\dblX$ be the involution interchanging $\Xsp_1$ and $\Xsp_2$ by the identity map.

Then the foliation $\Partition$ on each of the copies of $\Xsp$ gives a one-dimensional foliation $\dblP$ on $\dblX$.
Moreover, let $\dblV$ be the double of $\crsDomain$ as in Lemma~\ref{lm:double_of_open_set}.
Then $\dblV$ is open in $\bR^n$ and the cross section $\crs$ naturally extends to the cross section $\dblc:\dblV \to\dblX$ of $\dblP$ such that $\dblc|_{\crsDomain} = \crs$ and $\sigma\circ\dblc = \dblc\circ\xi$.

Since $\partial\dblX=\varnothing$, it follows from the boundary-less case of Theorem~\ref{th:trivialization_of_cross_section} that the saturation $\SAT{\dblc(\dblV)}$ is open in $\dblX$ and foliated homeomorphic with $\bR\times\dblV$.
That homeomorphism induces a homeomorphism of the open subset
\[ \SAT{\crs(\crsDomain)} = \SAT{\dblc(\dblV)} \cap \Xsp_1\]
of $\Xsp_1$ onto $\bR\times\crsDomain$.
Theorem~\ref{th:trivialization_of_cross_section} is completed.

\section{Acknowledgments}
The authors are sincerely grateful to Olena Karlova for pointing out to Katetov-Tong theorem about characterization of normal spaces which leads us to the proof of Theorem~\ref{th:trivialization_of_cross_section} using E.~Michael selection theorem.

\def\cprime{$'$}
\providecommand{\bysame}{\leavevmode\hbox to3em{\hrulefill}\thinspace}
\providecommand{\MR}{\relax\ifhmode\unskip\space\fi MR }
\providecommand{\MRhref}[2]{%
  \href{http://www.ams.org/mathscinet-getitem?mr=#1}{#2}
}
\providecommand{\href}[2]{#2}

\end{document}